 \documentclass[12pt]{article}
 \usepackage{amssymb}

\makeatletter
 
\@addtoreset{equation}{section}

\makeatother

\renewcommand\thefigure{\thesection.\@arabic\c@figure}

\renewcommand\thetable{\thesection.\@arabic\c@table}

\def\reff#1{(\ref{#1})}

\begin{document}

\def\E{{\Bbb E}}
\def\P{{\Bbb P}}
\def\R{{\Bbb R}}
\def\Z{{\Bbb Z}}
\def\V{{\Bbb V}}
\def\N{{\Bbb N}}
\def\X{{\cal X}}
\def\W{{\cal W}}
\def\G{{\cal G}}
\def\T{{\cal T}}
\def\I{{\cal I}}
\def\TT{\bar{{\cal T}}}
\def\II{\bar{{\cal I}}}
\def\C{{\C}}
\def\C{{\cal D}}
\def\n{{\bf n}}
\def\m{{\bf m}}
\def\b{{\bf b}}
\def\Var{{\hbox{Var}}}
\def\Cov{{\hbox{Cov}}}

\def\sqr{\vcenter{
         \hrule height.1mm
         \hbox{\vrule width.1mm height2.2mm\kern2.18mm\vrule width.1mm}
         \hrule height.1mm}}                  
\def\square{\ifmmode\sqr\else{$\sqr$}\fi}
\def\one{{\bf 1}\hskip-.5mm}
\def\limn{\lim_{N\to\infty}}
\def\given{\ \vert \ }
\def\ze{{\zeta}}
\def\be{{\beta}}
\def\la{{\lambda}}
\def\ga{{\gamma}}
\def\th{{\theta}}
\def\proof{\noindent{\bf Proof. }}
\def\A{{\bf A}}
\def\B{{\bf B}}
\def\C{{\bf C}}
\def\D{{\bf D}}
\def\MM{{\bf m}}
\def\w{\bar{w}}
\def\lnt{{\Lambda^N}}
\def\dlnt{\delta\Lambda^N_t}
\def\lno{\Lambda^N_0}
\def\dlno{\delta\Lambda^N_0}

\title{Flux fluctuations in the one dimensional
  nearest neighbors symmetric simple exclusion process}

\author{{\large A. De Masi}, {\it Universit\`a de L'Aquila}\\
{\large P. A. Ferrari}, {\it Universidade
de S\~{a}o Paulo}
}
\date{}
\maketitle
\noindent {\bf Abstract}
Let $J(t)$ be the the integrated flux of particles in the symmetric
simple exclusion process starting with the product invariant measure
$\nu_\rho$ with density $\rho$. We compute its rescaled asymptotic variance:
\[ \lim_{t\to\infty} t^{-1/2} \V J(t) \,=\, \sqrt{2/\pi}
(1-\rho)\rho 
\]
Furthermore we show that $t^{-1/4}J(t)$ converges weakly to a centered normal
random variable with this variance. From these results we compute the
asymptotic variance of a tagged particle in the nearest neighbor case and show
the corresponding central limit theorem.



\vskip 3mm

\paragraph{Results}
The nearest neighbors symmetric simple exclusion process describes the
evolution of particles sitting at the sites of $\Z$ evolving as follows. At
most one particle is allowed at each site. If there is a particle at a given
site, at rate one the particle chooses one of its nearest neighbor sites with
probability $1/2$ and attempts to jump to this site. The jump is effectively
realized if the destination site is empty; if not, the jump is suppressed. A
formal definition using Poisson processes is given below. The generator of the
process is given by
\begin{equation}
  \label{l1}
  Lf(\eta)\; =\; \frac12\sum_{x\in\Z}  [f(\eta^{x,x+1}) - f(\eta)]
\end{equation}
where $\eta^{x,x+1}(x) = \eta(x+1)$, $\eta^{x,x+1}(x+1) = \eta(x)$ and
$\eta^{x,x+1}(y) = \eta(y)$ for $y\neq x, x+1$.  For each $\rho\in[0,1]$ the
product measure $\nu_\rho$ with density $\rho$ is invariant for the process.

For an initial configuration $\eta$ let the integrated flux of particles
$J(t)= J^\eta(t)$ be the number of particles to the left of the origin at time
zero and to the right of it at time $t$ minus the number of particles to the
right of the origin at time $0$ and to the left of it at time $t$.

Fix $\rho\in (0,1)$ and let the initial configuration have law
$\nu_\rho$. Let
\[\V J(t) = \E^{\nu_\rho}J(t)^2=:\int d\nu_\rho(\eta) (J^\eta_t)^2.\] 
(Notice that $\E^{\nu_\rho}J(t)=0$.)

We prove the following asymptotics for the variance:
\begin{equation}
  \label{a}
  \lim_{t\to\infty} t^{-1/2} \V J(t) \,=\, \sqrt{2/\pi}
  (1-\rho)\rho \;:=\;\sigma^2
\end{equation}
We then prove the following central limit theorem for the integrated flux: 
\begin{equation}
  \label{17}
  t^{-1/4}J(t) \mbox{ converges weakly to } {\cal N}(0,\sigma^2) 
\end{equation}
where ${\cal N}(0,\sigma^2)$ is a centered normal random variable with
variance $\sigma^2$. 

Finally, let $X(t)$ be the position of a tagged particle interacting by
exclusion. We show that if the initial configuration is chosen with the
product measure $\nu_\rho$, then 
\begin{equation}
  \label{18}
 \lim_{t\to\infty} t^{-1/2} \E^{\nu_\rho}(X(t) -\rho^{-1} J(t))^2 = 0
\end{equation}
An immediate consequence of \reff{a}, \reff{17} and \reff{18} is that,
defining $ \V X(t) = \E^{\nu_\rho}(X(t))^2$, the
asymptotic variance of the tagged particle is 
\begin{equation}
  \label{29}
  \lim_{t\to\infty} t^{-1/2} \V X(t) \,=\, \sqrt {2/\pi} 
  {1-\rho\over \rho} \;:=\;\bar\sigma^2
\end{equation}
and the tagged particle satisfies a central limit theorem:
\begin{equation}
  \label{30}
  t^{-1/4}X(t) \mbox{ converges weakly to } {\cal N}(0,\bar\sigma^2) 
\end{equation}
The limits \reff{29} and \reff{30} were proven by Arratia (1983).

To prove the above results we use the stirring motion representation of the
symmetric exclusion process introduced by Harris (1972) and used by Arratia to
prove \reff{29} and \reff{30}.

\paragraph{The stirring process} The stirring process $z(i,t)\in \Z$, $i\in
\Z$, is defined as follows. At time $t=0$ put a (labeled) particle at each
site and define $z(i,0)=i$ for all $i\in \Z$. With each bond $(x,x+1)$, $x\in
\Z$ associate a Poisson process (clock) with parameter $1/2$. When the clock
rings at the bond $(x,x+1)$ the particles at those sites interchange their
positions. $z(i,t)$ is the position at time $t$ of the particle sitting at $i$
at time $0$. Given an initial configuration $\eta\in{\cal X}$, it is possible
to define the simple exclusion process $\eta_t$ in terms of the stirring
process by setting
\begin{equation}
  \label{1}
  \eta_t(x) = \one\{x\in \{z(i,t): \eta(i) =1\}\} 
\end{equation}

\paragraph{First proof of \reff{a}}
In terms of the stirring process, we define the following random variables.
\begin{equation}
  \label{4}
  K^+(t) = \sum_{i\le 0} \one\{z(i,t)>0\}\,;\qquad 
K^-(t) = \sum_{i> 0} \one\{z(i,t)\le 0\}
\end{equation}
where $\one\{\cdot\}$ is the characteristic function of the set $\{\cdot\}$.
The variable $K^+(t) $ is the number of stirring particles starting at the
left of the point $1/2$ and sitting at time $t$ at the right of $1/2$. The
variable $K^-(t) $ is the number of stirring particles starting at the right
of the point $1/2$ and sitting at time $t$ at the left of $1/2$. Since at all
times all sites are occupied by one stirring particle, each crossing of the
point $1/2$ from left to right involves a simultaneous crossing in the
opposite direction and viceversa. So $K^+(t)-K^-(t)$ is constant in $t$ and
since $K^+(0)=K^-(0)=0$, $K^+(t) = K^-(t) := K(t)$, for all $t\ge 0$. In the
stirring process the representation of $J(t)$ is given by
\begin{equation}
  \label{5}
  J(t) \;=\; \sum_{i\le 0} \one\{z(i,t)>0\}\eta(i)\,-\,
\sum_{i> 0} \one\{z(i,t)\le 0\}\eta(i)
\end{equation}
Let $i_1<i_2<\dots<i_{K(t)}\le 0$ be the random sites for which $z(i_k,t)>0$
  and $0<j_1<j_2<\dots<j_{K(t)}$ be the random sites for which
  $z(j_k,t)\le 0$. Define $B^+(k)= \eta(i_k)$ and $B^-(k)= \eta(j_k)$ and
  $A(k) = B^+(k)-B^-(k)$. Thus
    \begin{equation}
      \label{6}
      J(t) \;=\; \sum_{k=1}^{K(t)}A(k)
    \end{equation}
Assume $\eta$ is distributed according to the product measure $\nu_\rho$. Then
the variables $B^+(k)$, $B^-(k)$ and $K(t)$ are independent. Hence $A(k)$ are
iid independent of $K(t)$ with law
\begin{equation}
  \label{7}
  \P(A(k) = 1) = \P(A(k)=-1)= \rho(1-\rho)\;;\qquad \P(A(k) =0)=1-2\rho(1-\rho)
\end{equation}
Thus $\E A(k)=0$, $\E A(k)^2 = 2\rho(1-\rho)$ and by independence, using
\reff{6} we have
\begin{equation}
  \label{8}
  \E^{\nu_\rho} J(t)^2 = \E A(k)^2 \,\E K(t)
\end{equation}
To compute $\E K(t)$ write
\begin{eqnarray}
  \label{9}
  \E K(t) &=& \sum_{i\le 0} \P(z(i,t)>0)\;=\; \sum_{i\ge 0}
  \P(z(0,t)>i)\;=\; \E(z(0,t))^+.\nonumber  
\end{eqnarray}
But $z(0,t)$ is a symmetric random walk, thus, since $t^{-1}\E z(0,t)^2$ is
uniformly integrable,
\begin{equation}
  \label{10}
  \lim_{t\to\infty} t^{-1/2}\E K(t)\;=\;  {1\over \sqrt
  {2\pi}}
\end{equation}
Thus, using \reff{8} we obtain \reff{a}. \square

\paragraph{Second proof of \reff{a}}
 From the definition we have 
\begin{equation}
  \label{11}
  J(t) - \int_0^t {1\over 2}(\eta_s(0) - \eta_s(1)) := M(t)
\end{equation}
where $M(t)$ is a martingale with variance
\begin{equation}
  \label{12}
  \E^{\nu_\rho} M(t)^2 = t\rho(1-\rho)
\end{equation}
As in De Masi {\it et al} (1985,1989), from the time invariance of
$\nu_\rho$
and the fact that $  J(t)$ is an anti-symmetric random variable, it
follows that
\begin{equation}
  \label{13}
  \E^{\nu_\rho} J(t)^2 = t\rho(1-\rho)- {1\over 2}\int_0^t ds (t-s)
 \int \nu_\rho(d\eta) (\eta(0) - \eta(1))\E(\eta^\eta_s(0) -
\eta^\eta_s(1))
\end{equation}
where $\eta^\eta_s$ is the exclusion process with initial configuration
$\eta$.  {From} the reversibility and the translation invariance of
$\nu_\rho$,
\begin{equation}
  \label{14}
  \int \nu_\rho(d\eta) (\eta(0) - \eta(1))\E(\eta^\eta_s(0) -
  \eta^\eta_s(1))=
2 \int \nu_\rho(d\eta) \Bigl(\eta(0) \E\eta^\eta_s(0) -
  \eta(0) \E\eta^\eta_s(1)\Bigr)
\end{equation}
Calling $L$ the generator of the process we have that
\begin{equation}
L\eta(0)=\frac{1}{2}[\eta(1)-\eta(0)]+\frac{1}{2}[\eta(-1)-\eta(0)]
\end{equation}
Therefore, using once more translation invariance
\begin{equation}
  \label{14a}
2 \int \nu_\rho(d\eta) \Bigl(\eta(0) \E\eta^\eta_s(0) -
  \eta(0) \E\eta^\eta_s(1)\Bigr)
= -2 {d\over ds} \Bigl(\int \nu_\rho(d\eta) \,\eta(0)
  [\E\eta^\eta_s(0)-\rho]\Bigr) 
\end{equation}
We use  \reff{14} and \reff{14a} in the second term on the right hand side
of
\reff{13} then,  integrating by parts, we get
  \begin{equation}
    \label{15}
    \int_0^t (t-s) {d\over ds} \Bigl(\int \nu_\rho(d\eta) \,\eta(0)
  \E\eta^\eta_s(0)-\rho^2\Bigr)\;
= \;-t \rho(1-\rho) + \int_0^t \int \nu_\rho(d\eta) \,(\eta(0)-\rho)
  \E\eta^\eta_s(0)
  \end{equation}
{From} \reff{13} and \reff{15} we finally get
\begin{eqnarray}
  \label{16}
  \E^{\nu_\rho} J(t)^2 &=& \int_0^t \int \nu_\rho(d\eta) \,(\eta(0)-\rho)
  \E\eta^\eta_s(0)\\
&=& \rho(1-\rho) R_t(0)
\end{eqnarray}
where $R_t(0)$ is the expected amount of time spent at the origin up to
time
$t$ for a continuous time symmetric random walk starting at zero. Finally,
\begin{equation}
  \label{37}
 \lim_{t\to\infty} t^{-1/2} R_t(0) = \sqrt{2/\pi}.\qquad\square
\end{equation}

\paragraph{Proof of \reff{17}}
To show \reff{17} from \reff{6} it is enough to show that 
\begin{equation}
  \label{20}
  C(t) := t^{-1/4} \Bigl(\sum_{k=1}^{K(t)}A(k) -
  \sum_{k=1}^{t^{1/2}/\sqrt{2\pi}}A(k)\Bigr) \rightarrow 0 \mbox{ as }
  {t\to\infty} 
\end{equation}
in measure. Using Chebishev inequality we have, for any $c>0$,
\begin{equation}
  \label{21}
  \P(C(t)>c) \le {\E A(k)^2\over c^2} \E\Bigl|{K(t)\over t^{1/2}} - {1\over
  \sqrt{2\pi}}\Bigr| \rightarrow 0 \mbox{ as }
  {t\to\infty}
\end{equation}
The limit goes to zero because $K(t)$ is the sum of negatively correlated
0--1 random variables and so $\V K(t) \le \E K(t) \sim \sqrt t$ (Arratia
(1983)) and by Schwarz inequality.

\paragraph{Proof of \reff{18}}
We use a lattice version of a result of D\"urr, Goldstein and Lebowitz (1985)
for an infinite ideal gas of point particles on $\R$. Suppose that the initial
configuration $\eta$ is distributed according to the invariant measure
$\nu_\rho$. Fix $t\ge 0$. For $k\ge 0$ let $Y_k(t)$ be the position of the
$k$th particle of $\eta_t$ to the right of $1/2$, with $Y_0(t)\le 0$. For $k<
0$ let $Y_k(t)$ be the position of the $-(k+1)$th particle of $\eta_t$ to the
left of $1/2$. (When time goes on the particles change these labels.) It is
easy to see that at time $t$ th
e tagged particle (which at time $t=0$ is
labeled 0) is the $J(t)$th particle, that is:
\begin{equation}
  \label{22}
  X(t) = Y_{J(t)}(t)
\end{equation}
By the ergodicity (under translations) and stationarity of $\nu_\rho$ we have
that 
\begin{equation}
  \label{23}
  \lim_{n\to\infty} n^{-1} Y_n(t) = \rho^{-1},\qquad \P^{\nu_\rho}-\mbox{almost
  surely}. 
\end{equation}
One can then prove (as in Lemma 2.8 of D\"urr, Goldstein and Lebowitz (1985))
that 
\begin{equation}
  \label{24}
  \lim_{t\to\infty} t^{1/2} \E^{\nu_\rho}( Y_{J(t)}(t) - \rho^{-1} J(t))^2 = 0
\end{equation}

\paragraph{Acknowledgments.} 
We thank discussions with Errico Presutti, Shelly Goldstein and David Wick.
PAF thanks support from CNPq.

\paragraph{Remark} 
This work was written when the authors visited Rutgers University in 1985 and
was kept unpublished for more than 15 years. We decided to publish it now for
three reasons. The first proof of \reff{a} is an application of Arratia's
method, but it is not written anywhere; in fact, it is easier first to compute
the variance of the flux and then, as a corollary, the variance of the tagged
particle than vice versa. The second proof of \reff{a} is the unique
application we know of the method of De Masi et al. (1985--1989) that works
for a subdiffussive process. Finally, the flux in the simple exclusion process
is isomorphic to a 1+1 dimensional interface. The role of the entropic
repulsion when this interface interacts by exclusion with a wall has been
studied by Dunlop, Ferrari and Fontes (2001), who compare the asymptotic
variance of the flux for the process starting with the deterministic
configuration $\dots101010\dots$ with the stationary process studied here.

Presumably our result can be obtained using the fact that the asymptotic
behavior of the current can be deduced from the hydrodynamic behavior of the
symmetric simple exclusion and the asymptotics of the density fluctuation
field at equilibrium. This technique has been introduced by Rost and Vares
(1985) and applied to the zero range process by Landim, Olla and Volchan
(1997--1998) and Landim and Volcham (2000).  We are not aware of any
application of this argument to our case and after consulting Landim and Olla
it seems that their results do not cover, at least automatically, ours. We
thank an anonymous referee, Errico Presutti, Claudio Landim and Stefano Olla
for pointing out this possibility.

\vskip 3mm
\parindent 0pt
Anna De Masi \\
Dipartimento di Matematica Pura ed Applicata\\
Universita' di L'Aquila
I- 67100 - L'Aquila - Italy - \\
email: {\tt demasi@univaq.it}

\vskip 4mm
Pablo A. Ferrari \\
IME USP, Caixa Postal 66281, 05315-970 - S\~{a}o Paulo, SP - BRAZIL - \\
email: {\tt pablo@ime.usp.br}\\


\begin{thebibliography}{999}
  
\bibitem{a} Arratia, Richard (1983) The motion of a tagged particle in the
  simple symmetric exclusion system on ${Z}$.  {\sl Ann. Probab. \bf 11}, no.
  2, 362--373.

\bibitem{dfgw1} De Masi, A.; Ferrari, P. A.; Goldstein, S.; Wick, W. D. (1985)
  Invariance principle for reversible Markov processes with application to
  diffusion in the percolation regime. {\sl Particle systems, random media and
    large deviations} (Brunswick, Maine, 1984), 71--85, {\sl Contemp. Math.},
  41, Amer.  Math. Soc., Providence, R.I.

\bibitem{dfgw1}De Masi, A.; Ferrari, P. A.; Goldstein, S.; Wick, W. D. 
  (1989) An
  invariance principle for reversible Markov processes.  Applications to
  random motions in random environments. {\sl J. Statist. Phys. \bf 55},
  no. 3-4, 787--855.
 

\bibitem{dff} Dunlop, F. M. ; Ferrari, P. A.; Fontes, L. R. G. (2001)
A dynamic one-dimensional interface interacting with a wall.
To appear in J. Stat. Phys. 

  
\bibitem{dgl} D\"urr, Detlef; Goldstein, Sheldon; Lebowitz, Joel L.  (1985)
  Asymptotics of particle trajectories in infinite one-dimensional systems
  with collisions. {\sl Comm. Pure Appl. Math. \bf 38},

\bibitem{h} Harris, T. E.  (1972)
Nearest-neighbor Markov interaction processes on multidimensional lattices. 
{\sl Advances in Math. \bf 9}, 66--89. 


\bibitem{lvo1} Landim, C.; Olla, S.; Volchan, S. B.  (1997) Driven tracer
  particle and Einstein relation in one-dimensional symmetric simple exclusion
  process.  {\sl Resenhas \bf 3},2:173--209.
  
\bibitem{lvo2} Landim, C.; Olla, S.; Volchan, S. B.  (1998) Driven tracer
  particle in one-dimensional symmetric simple exclusion. {\sl Comm. Math.
    Phys.  \bf 192} 2:287--307.

\bibitem{lv} Landim, C.; Volchan, S. B.  (2000) Equilibrium
  fluctuations for a driven tracer particle dynamics. {\sl Stochastic Process.
  Appl. \bf 85} 1:139--158.


\bibitem{LOV} Landim, C.; Olla, S.; Volchan, S. B.(1997)
 Driven tracer particle
 and Einstein relation in one-dimensional
 symmetric simple exclusion process.
 {\sl Resenhas \bf 3}, 173--209

\bibitem{RV} Rost, H.. Vares M.E.  (1985)
Hydrodynamics of a one-dimensional nearest neighbor model.
Particle systems, random media and large deviations . {\sl
Contemp. Math. \bf {41}}.


\end{thebibliography}
\end{document}